\documentstyle[12pt,amsthm,amssymb,amscd,fullpage]{amsart}
\input BoxedEPS.tex
\SetRokickiEPSFSpecial
\HideDisplacementBoxes
\def\draft{y}
\def\printname#1{
	\if\draft n
		\smash{\makebox[0pt]{\hspace{-0.5in}
			\raisebox{8pt}{\tt\tiny #1}}}
	\fi
}
\def\lbl#1{\label{#1}\printname{#1}}
\theoremstyle{plain}

\newtheorem{thm}{Theorem}
\newtheorem{prop}{Proposition}[section]
\newtheorem{lemma}[prop]{Lemma}
\newtheorem{cor}[prop]{Corollary}

\theoremstyle{definition}

\theoremstyle{remark}

\catcode`\@=12

\def\a{\alpha}
\def\b{\beta}
\def\S{\Sigma}

\def\th{\theta}
\def\Th{\Theta}
\def\s{\sigma}

\def\n{\nabla}
\def\l{\lambda}
\def\L{\Lambda}
\def\m{\mu}
\def\O{\Omega}
\def\Ga{\Gamma}
\def\gg{\gamma}
\def\d{\delta}
\def\D{\Delta}

\def\sub{\subseteq}
\def\pa{\partial}
\def\sh{\sharp}

\def\TT{{\cal T}}

\def\LL{{\cal L}}
\def\SS{{\cal S}} 
\def\I{{\cal I}}

\def\Z{{\Bbb Z}}
\def\ZZ{\bold Z}
\def\Q{{\Bbb Q}}
\def\R{{\Bbb R}}

\def\ld{\ldots}

\def\im{\operatorname{Im}}
\def\int{\operatorname{int}}

\def\iso{\cong}

\def\zz{\Z [z]}
\def\zs{\Z [s,s^{-1}]}
\def\zss{\Z [s,s^{-1}]_{\Sigma '}}
\def\zt{\Z [t,t^{-1}]}
\def\zts{\zt_{\S}}
\def\zti{\Z [t_1 ,\cdots ,t_m ,t_1^{-1},\cdots ,t_m^{-1}]}
\def\zf{\Z [F]}
\def\zu{\Z [[u]]}
\def\zui{\Z [[u_1 ,\ldots ,u_m ]]}
\def\zvi{\Z [[v_1 ,\ldots ,v_m ]]}
\def\X{[X_1 ,\ld ,X_m ]}
\def\nl{\nabla_L (z)}
\def\mb{\bar\mu}
\def\ti{\tilde}
\def\mb{\bar\mu}

\def\tith{\ti\th}
\def\ap{\approx}
\def\ba{\bar A}

\begin{document}

\title[A factorization of the Conway polynomial]{A factorization of the
Conway polynomial}
\author{Jerome Levine}
\thanks{Partially supported by NSF grant 
       DMS-96-26639. This and other preprints of the author can be
obtained  on the WEB at the address 
 {\tt 
http:\linebreak[0]//\linebreak[0]www.\linebreak[0]math.\linebreak
[0]brandeis.\linebreak[0]edu/\linebreak[0]Faculty/\linebreak[0]jlevine/.}}
\address{Department of Mathematics\\
        Brandeis University\\
        Waltham, MA 02254-9110, U.S.A. }
\email{levine@binah.cc.brandeis.edu}

\date
{ \today } 
\maketitle
\section{Introduction} It is tempting to conjecture that there is some
interesting relationship between the Conway polynomial $\n_L (z)$  of a
link $L$ and $\n_K (z)$, where $K$ is a knot obtained by banding
together the components of $L$. Obviously they cannot be equal since
only terms of even or odd degree appear in $\n_L
(z)$, according to whether $L$ has an
odd or even number of components.
Moreover there are many ways of choosing
bands and one can easily see that the
variety of knots one obtains can have
very different polynomials. Nevertheless
we will demonstrate that $\n_L (z)$ and
$\n_K (z)$ have a very precise
relationship in the form of a
factorization:
$$ \n_L (z)=\n_K (z)\Ga (z) $$
where $\Ga (z)$ is a power series in $z$ which depends on the choice
of bands. More precisely the choice of bands can be viewed as the choice
of a {\em string link }representation of $L$ and the coefficients of
$\Ga (z)$ are given by explicit formulae in terms of the Milnor
$\mb$-invariants of this string link. Thus it actually only depends on
the {\em I-equivalence class } of  $L$ (and the bands). From this point
of view we see that the indeterminacy of the band choice is compensated
by the notorious and familiar indeterminacy of the $\mb$-invariants.
However there is no indeterminacy in the first non-zero coefficient of
$\n_L (z)$ and our results give a general formula for this  coefficient
in terms of the $\mb$-invariants of $L$, generalizing the special
cases obtained previously by Cochran, Hoste and myself.

We will also obtain an analogous factorization of the {\em
multivariable} Alexander polynomial, depending on a choice of string
link $S$ representing $L$:
$$\D_L (t_1 ,\cdots ,t_m )=\th (t_1 ,\cdots ,t_m )\Ga (t_1 ,\cdots
,t_m ) $$
in which $\th (t_1 ,\cdots ,t_m )$ and $\Ga (t_1 ,\cdots
,t_m )$ are rational power series in $\{ t_i -1\}$ with the
properties that $\Ga (t_1 ,\cdots ,t_m )$ is given by an explicit
formula in terms of the $\bar\mu$-invariants of $S$ and $\th
(1,\cdots ,1)=1$. However in this factorization we cannot tell if
$\th (t_1 ,\cdots ,t_m )$ is a polynomial and, in any case, it
lacks the geometric interpretation of the corresponding factor in
the one-variable case. One interesting consequence though is a
formula for the lowest degree terms in the Taylor expansion of $\D_L
(t_1 ,\cdots ,t_m )$ about $(1,\cdots ,1)$ in terms of the
$\bar\mu$-invariants of $S$. These results are at least implicit
in Traldi \cite{Tra} and are closely related to work of Rozansky
\cite{Roz}. It is suggested by recent work of Habegger and Masbaum
\cite{HM} that such a factorization occurs for the general class of
finite-type invariants of string links.

Finally I would like to thank Stavros Garoufalidis for helpful
discussions.

\subsection{Statement of results} Suppose $L$ is an $m$-component
oriented link in $\R^3$. It is obvious that $L$ can be obtained by
{\em closing } a string link $S$-- see Section \ref{sec.str} for the
definitions. In fact there are many choices of different string links whose
closure give $L$ if $m>1$. Now for any string link $S$ we will define, in
Section
\ref{sec.ks}, another form of closure which will produce a {\em knot} which
we denote $K_S$. We will call $K_S$ the {\em knot closure} of $S$. In
fact $K_S$ will be a {\em band sum} of the components of $L_S$ (the usual
link closure of $S$). It is not hard to see that, for any knot $K$ obtained by
band-summing of the components of a link $L$, there is some string link
$S$ such that $L=L_S$ and $K=K_S$. 

In \cite{Mil} Milnor defined, for any oriented link $L$, an array of
integer-valued invariants $\{\mb_{i_1 ,\ldots ,i_k}(L)\}$. The study of
these invariants has always been hampered by a complicated self-referencing
indeterminacy in their definition. The simplest case  when they are
well-defined is formulated in the following recursive way: $\{\mb_{i_1
,\ldots ,i_k}(L)\}$ is well-defined if $\{\mb_{j_1 ,\ldots ,j_r}(L)\}$ is
well-defined and {\em zero} for any proper, order-preserving subset $\{ j_1
,\ldots ,j_r\}$ of $\{ i_1 ,\ldots ,i_k\}$. So, for example, the
$\mb$-invariants of order $k$ are well-defined if all those of order
less than $k$ are well-defined and zero. The situation was improved by the
introduction of string links. For a string link $S$, $\mb_{i_1 ,\ldots
,i_k}(S)$ is well-defined with {\em no indeterminacy}-- the definition is
recalled in Equation \eqref{eq.mb}. Thus the indeterminacy of the
$\mb$-invariants of
$L$ reflects the indeterminacy in the choice of a string link whose
closure is $L$.

Our main result will be:
\begin{thm}\lbl{th.slnk} Let $S$ be a string link, with closure $L_S$
and knot closure $K_S$. Then we have the following factorization of the
Conway polynomial of $L_S$:
$$ \n_{L_S}(z)=\n_{K_S}(z)\Ga_S (z) $$
where $\Ga_S (z)$ is a power series given by the formula:
$$\Ga_S (z)=(u+1)^{e/2}\det (\l_{ij}(u)) $$
with $z=u/\sqrt{u+1}$, $e=\cases 0\quad\text{if $m$ is odd}\\
1\quad\text{if $m$ is even}\endcases$ and:
  $$\l_{ij}(u)=\sum_{r=0}^{\infty}(\sum_{i_1 ,\ldots ,i_r}\mb_{i_1 ,\ldots
,i_r ,j,i}(S))\ u^{r+1} $$
\end{thm}

Note that $\Ga_S (z)$ is  a {\em rational} function, i.e. it is a quotient
of polynomials with integer coefficients. 

\begin{cor}\lbl{cor.lnk} If the $\mb$-invariants of $L$ vanish for order
less than $k$, and so the $\mb$-invariants of order $k$ are well-defined
(i.e. have no indeterminacy), then $\n_L (z)$ is divisible by
$z^{(k-1)(m-1)}$ and the coefficient of $z^{(k-1)(m-1)}$ is $\det
(a_{ij})$ where 
$$a_{ij}=\sum_{i_1 ,\ldots ,i_{k-2}}\mb_{i_1 ,\ldots ,i_{k-2},j,i}\ (L)
$$
\end{cor}
This generalizes results of Hoste \cite{Hos}, Cochran \cite{Coc} (also
see \cite[Corollary 6.3]{Tra}) and myself
\cite{Lev}.

\begin{cor}\lbl{cor.conc}
The first non-vanishing term (and its degree) of the Conway polynomial
is an I-equivalence invariant.
\end{cor}
This was first proved by Cochran \cite{Coc}. {\em I-equivalence } is
the relation generated by concordance and connected sum of a component
with a local knot. A {\em local knot } is a knot lying in a ball
disjoint from the link.

To state the next results we use the following terminology. A polynomial or
rational function $f(z)$ will be called {\em norm-like} if it can be
written in the form $f(z)=h(t)h(t^{-1})$ for some rational function
$h(t)$ ($z=t-t^{-1}$). If $f(z)$ is norm-like it must contain only
even powers of $z$ and $f(0)$ must be square. But, for example, $z^2
+a$ is norm-like if and only if $a=4$. Note that norm-like is not the
same as being of the form $g(z)g(-z)$.

\begin{cor}\lbl{cor.con} If $L$ and $L'$ are concordant links, then their Conway
polynomials are related in the following way:
$$\n_L (z)f(z)=\n_{L'}(z)g(z) $$
where $f,g$ are norm-like polynomials satisfying $f(0)=1=g(0)$.
\end{cor}
This was first proved in \cite{Coc}; the analogous result for the
multi-variable Alexander polynomial was first  proved in 
\cite{Kaw} and \cite{Nak}.

We will also prove the following result about $\Ga_S (z)$.
\begin{prop}\lbl{prop.gam}
Suppose that all the linking numbers of $L_S$ are zero. Then:
\begin{enumerate}
\item $\Ga_S (z)$ has only even (resp., odd) powers of $z$ if $m$ is
odd (resp., even).
\item $\Ga_S (z)$ is divisible by $z^{2(m-1)}$. If $m$ is odd, then either
$\Ga_S (z)$ is divisible by $z^{2m}$ or it is norm-like.
\end{enumerate}
\end{prop}
\begin{cor}
If all the linking numbers of $L_S$ are zero then the value of \ 
$z^{-2(m-1)}\Ga_S (z)$ at $z=0$ is square.
\end{cor}
This was recently proved in \cite{Lev}.

For the {\em multivariable} Alexander polynomial we will prove:
\begin{thm}\lbl{th.mual} There exists a rational power series $\th
(v_1 ,\cdots ,v_m )$ with constant term $\pm 1$, such that
$$\D_{L_S}(v_1 ,\cdots ,v_m )=\Phi_S (v_1 ,\cdots ,v_m )\th (v_1
,\cdots ,v_m ) $$
where 
$$\Phi_S (v_1 ,\cdots ,v_m )=\det (v_i (\sum_{i_1 ,\cdots
,i_r}\mb_{i_1 ,\ldots ,i_r ,j,i}(S)\ v_{i_1}\cdots v_{i_r})-(\tau_i
-1)\delta_{ij}) $$
\end{thm}
\begin{cor} If the $\mb$-invariants of $L$ vanish for order
less than $k$ (assume $k\ge 2$), and so the $\mb$-invariants of order $k$
are well-defined (i.e. have no indeterminacy), then $\D_L(v_1 ,\cdots ,v_m )$
has no terms of degree $< (k-1)(m-1)$ and the homogeneous part of $\D_L(v_1
,\cdots ,v_m )$ of degree $(k-1)(m-1)$ is $\det (a_{ij})$, where:
$$a_{ij}= v_i (\sum_{i_1 ,\cdots
,i_{k-2}}\mb_{i_1 ,\ldots ,i_{k-2} ,j,i}(L)\ v_{i_1}\cdots v_{i_{k-2}}$$
\end{cor}
\subsection{Outline of proof} For the Conway polynomial we will produce
factorizations
$\n_L (z)=\n'(z)\n''(z)$ from two different points of view. In Sections
\ref{sec.Hom} and \ref{sec.Alex} we will examine the homology of the
infinite cyclic cover of a string link and, in Theorem \ref{th.Phi}, produce
a factorization in which $\n'(0)=1$ and $\n''(z)=\Ga_S (z)$. In
Section \ref{sec.ks} we will introduce the Seifert matrix of a string link
and obtain another factorization, in equation \eqref{eq.fact}, in which
$\n'(z)=\n_{K_S}(z)$ and $\n'' (z)$ is defined from a Seifert matrix.
The proof is completed by showing that these two factorizations are actually
the same. This argument is carried out in Section \ref{sec.Seif}.

For Theorem \ref{th.mual} we only need carry out the homological argument.
This is explained in Sections \ref{sec.unab} and \ref{sec.mual}.

In Section \ref{sec.prf} we prove Proposition \ref{prop.gam} and in Section
\ref{sec.ex} we illustrate the factorization of $\n_L (z)$ given in Theorem
\ref{th.slnk} by a very simple example.

\vspace{.2in}
\section{Preliminaries}
\subsection{String links}\lbl{sec.str}
We review the definition of the $\mb$-invariants of a string link. 

Recall that a string link $S$ is an ordered collection of disjoint
oriented, properly imbedded arcs $S_1 ,\ldots ,S_m$ in $I\times
D^2$. It will be convenient for our purposes to orient the strings in the
following manner. The odd numbered strings $S_i$ will be directed from $(0,p_i
)$ to $(1, p_i )$ and the even-numbered strings from $(1, p_i )$  to
$(0,p_i )$, where
$p_1 ,\ldots ,p_m$ are prescribed distinct points in $D^2$. The {\em closure }
$L_S$ of $S$ is the oriented link in $S^3$ whose components $L_i
=S_i\cup A_i$, where $A_i$ are prescribed disjoint arcs in
$\overline{S^3 -I\times D^2}$ connecting $(1,p_i )$ and $(0,p_i )$,
which meet $I\times D^2$ only at their endpoints. See Figure \ref{fig.1}.

\begin{figure}[hbt]
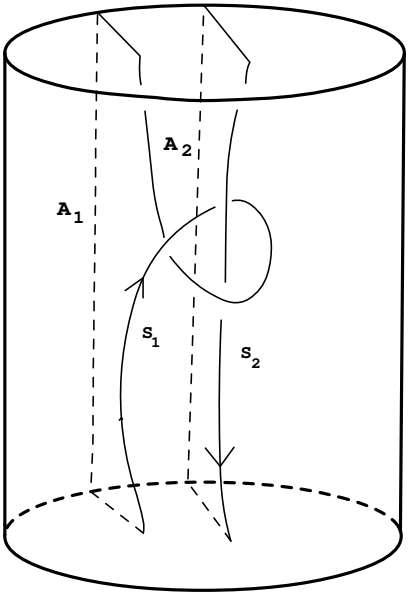

\begin{center}{\BoxedEPSF{Fig1 scaled 600}}\caption{Closure of a string
link}\lbl{fig.1}
\end{center}
\end{figure}
If $X=(I\times D^2 )-S$ and $\pi =\pi_1 (X)$, then $H_1 (\pi )$
is free abelian of rank $m$ generated by meridian elements,  and $H_2
(\pi )=0$. The latter follows since we can choose a Wirtinger
presentation for $\pi$ of deficiency $m$. We may choose {\em canonical
} meridian elements $\mu_1 ,\ldots ,\mu_m \in\pi$ represented by curves
lying in $0\times D^2$ (assuming the base point lies in $0\times S^1$) so that
the linking number of $\mu_i$ and $L_i$ is $+1$. If $F$ is the free group on
generators
$x_1 ,\ldots ,x_m$, then the map
$F\to\pi$  defined by $x_i\to\mu_i$ induces isomorphisms $F/F_q\iso\pi
/\pi_q$, for all $q$, and $\tilde F\to\tilde\pi$, where, for any
group $G$, $\{ G_q\}$ denotes the lower central series of $G$-- $G_1 =G$
and $G_{q+1}=[G.G_q ]$-- and $\tilde G=\varprojlim G_i$, the nilpotent
completion. This follows from Stallings theorem \cite{Sta}. We can also
choose canonical {\em longitudes }
$\l_i\in\pi$ with representative curves parallel to the components of
$S$ and closed up in the boundary of $I\times D^2$ in a prescribed
way, so that the {\em total } linking number of $\l_i$  is
zero. The total linking number of any oriented closed curve $\a$ in
$X$ is defined to be the sum of the linking numbers of $\a$ with
all the components of $L_S$. This choice of $\l_i$ differs from
the traditional choice in which the linking number of $\l_i$ and
$L_i$ is required to be zero, but our choice will be more convenient
for formulae. If the link is {\em algebraically split }, i.e. all the
linking numbers of the components of $L_S$ are zero, then, of
course, both choices are the same.

\subsection{Magnus expansion and the $\mb$-invariants}
In \cite{Mag} Magnus defines an imbedding $\th :F\to \zui$, where $\zui$
is the power-series ring  in $m$ non-commuting
variables $u_1 ,\ld ,u_m$, by defining $\th (x_i )=1+u_i$ and
extending this to a group homomorphism into the multiplicative group
$1+\I$. $\I$ is the ideal in $\zui$ consisting of all series
with zero constant term. Since
$\th (F_q )\sub 1+\I^q
$,
$\th$ extends to define an imbedding $\tilde F\to \zui$. For a string
link
$S$ this induces an imbedding $\th_S :\tilde\pi\to \zui$. We can then
extend
$\th$ and
$\th_S$ to ring homomorphisms $\Z F\to \zui$ and $\Z\pi\to \zui$. Now,
following the original concept in Milnor \cite{Mil}, but adapting to the
string link context (see e.g. \cite{Le2}), we define the
$\mb$-invariants of
$S$ by the formula
\begin{equation}\lbl{eq.mb}
\th_S (\l_i )=1+\sum_{i_1 ,\ld ,i_r}\mb_{i_1 ,\ld ,i_r ,i}\
(S)u_{i_1}\cdots u_{i_r}
\end{equation}
This definition will differ, in an easily formulated way (which we
omit), from the traditional definition, because we have used a
different choice of $\l_i$. If
$L_S$ is the closure of
$S$, then it is clear that
$\mb_{i_1 ,\ld ,i_r}\ (S)\equiv\mb_{i_1,\ld ,i_r}\ (L_S)$, modulo
indeterminacy, where  $\mb_{i_1,\ld ,i_r}\ (L_S)$ denotes the classical
$\mb$-invariants of $L_S$ defined by Milnor in \cite{Mil} (again with the
reservation mentioned above).
\section{Homological invariants of string links}\lbl{sec.Hom}
\subsection{The longitudinal matrix}\lbl{sec.long}
Consider now the infinite cyclic cover $p:\tilde X\to X$ defined by
the epimorphism $\pi\to\ZZ$ sending $\mu_i\to t$, where $\ZZ$ is
the infinite cyclic multiplicative group with generator $t$. We
consider the relative homology group $H_1 (\ti X,\ti\ast )$, where
$\ast$ is the base-point of $X$ and $\ti\ast =p^{-1}(\ast )$.
This is a module over the Laurent polynomial ring $\zt$. Consider
the multiplicative subset $\S\sub\zt$ consisting of all $f(t)$
such that $f(1)=1$, and the {\em localization } $\zts$ of $\zt$
consisting of all quotients of elements of $\zt$ by elements of
$\S$. We could, alternatively, consider the {\em completion } of
$\zt$. Let $\zu$ denote the ring of power series in the variable
$u$. Then we can define a ring homomorphism $\ti\th :\zt\to\zu$ by
$\ti\th (t)=1+u$. Since $\ti\th (\S )\sub\text{ units of }\zu$,
$\ti\th$ extends to an imbedding $\zts\to\zu$ which we still
denote by $\tith$.
\begin{lemma}\lbl{lem.free} $H_1 (\ti X,\ti\ast )_\S$ is a free
$\zts$-module with basis $\ti\mu_1 ,\ld ,\ti\mu_m$, where
$\ti\mu_i$ is represented by a lift of $\mu_i$ to a curve in $\ti
X$ starting at some prescribed base-point $\hat\ast\in\ti\ast$.
\end{lemma}

For a proof see \cite{Le3}. This lemma will also follow from the argument
given below in Section \ref{sec.seif}.

Now suppose we define $\ti\l_i\in H_1 (\ti X,\ti\ast )$ to be the
class represented by the lift of $\l_i$ to a path beginning at
$\hat\ast$. Then we can write, in $H_1 (\ti X,\ti\ast )_{\S}$:
\begin{equation}\lbl{eq.cij}
\ti\l_i =\sum_{j=1}^m c_{ij}^S\ti\mu_j ,\quad 1\le i\le m
\end{equation}
We will show how the $c_{ij}^S$ are determined by the
$\mb$-invariants of $S$.
\begin{lemma}\lbl{lem.cij} If we write
\begin{equation}\lbl{eq.thcij}
\ti\th (c_{ij}^S)=\sum_{k=0}^{\infty}c_{ijk}(S)u^k
\end{equation}
then
\begin{equation}\lbl{eq.cijk}
c_{ijk}(S)=\sum_{i_1 ,\ld ,i_k}\mb_{i_1 ,\ld ,i_k ,j,i}(S)
\end{equation}
\end{lemma}
\begin{pf} The relation between the image of  $\l_i$
under the mapping $\pi\to\ti\pi
\iso\ti F$ and $\ti\l_i\in H_1 (\ti X,\ti\ast )_{\S}\iso\zts
[\ti\mu_1 ,\ld ,\ti\mu_m ]$ is well-understood, on the purely
algebraic level, to be given by the Fox free differential calculus as
follows.  If $g\in F$, then consider $g-1\in IF\sub\zf$, where
$IF$ is the augmentation ideal in $\zf$. $IF$, regarded as a
left $\zf$-module, is freely generated by the elements $x_1 -1,\ld
, x_m -1$ and so we can write $g-1=\sum_i a_i (x_i -1)$. Let $\eta
:\zf\to\zt$ be defined by $\eta (x_i )=t$. Let $\ti g$ denote the 
element $\sum_i\eta (a_i )X_i\in\zt [X_1
,\ld ,X_m ]$, the free $\zt$-module with basis $X_1 ,\ldots ,X_m$. Then
$\d (g)=\ti g$ defines a function $\d :F\to\zt [X_1 ,\ld ,X_m ]$,
which is additive and has the property $\d (gh)=\eta (g)\d (h)+\d
(g)$. Also note that $\d (F_{q+1})\sub I^q\zt [X_1 ,\ld ,X_m ]$, where
$I$ is the augmentation ideal of $\zt$. It follows that $\d$ induces
a function $F/F_{q+1}\to (\zt /I^q )\X$.

Now if $\l_{iq}$ is the reduction of $\l_i$ into $\pi
/\pi_{q+1}\iso F/F_{q+1}$, then $\d (\l_{iq})$ is the reduction of
$\ti\l_i$ into $H_1 (\ti X,\ti\ast )_{\S}/I^q H_1 (\ti X,\ti\ast
)_{\S}\
\iso (\zt /I^q )\X$, where $\ti\mu_i\to X_i$.

Let us now write 
\begin{equation}\lbl{eq.liq}
\l_{iq}-1=\sum _j a_{ij}(x_j -1)
\end{equation}
 where the
$a_{ij}$ are well-defined mod $(IF)^q$. Then, by the discussion above,
we have 
$$\ti\l_i\equiv\sum_j\eta (a_{ij})\ti\mu_j\ \mod I^q $$
and so $c_{ij}^S\equiv\eta (a_{ij}) \mod \I^q$. From this we have 
\begin{equation}\lbl{eq.eta}
\ti\th (c_{ij}^S)\equiv\ti\eta (\th (a_{ij}))\ \mod u^q
\end{equation}
 where $\ti\eta
:\zui\to\zu$ is defined by $\ti\eta (u_i )=u$.

On the other hand if we apply $\th_S$ to equation \eqref{eq.liq} we
get
$$\th_S (\l_{iq})\equiv 1+\sum_j\th (a_{ij})u_j\ \mod I^{q+1}$$
Comparing this with equation \eqref{eq.mb} we have
$$\sum_j\th (a_{ij})u_j\equiv\sum_{i_1 ,\ld ,i_r}\mb_{i_1 ,\ld ,i_r ,i}
(S)u_{i_1}\cdots u_{i_r}\ \mod\I^{q+1} $$
and from this we see that
\begin{equation}\lbl{eq.aij}
\th (a_{ij})\equiv\sum_{i_1 ,\ld ,i_r}\mb_{i_1 ,\ld ,i_r
,j,i}(S)u_{i_1}\cdots u_{i_r}\ \mod \I^q
\end{equation}
Finally we combine equations \eqref{eq.eta} and \eqref{eq.aij} to get
\begin{equation}\lbl{eq.mij}
\ti\th (c_{ij}^S)\equiv\sum_{i_1 ,\ld ,i_r}\mb_{i_1
,\ld ,i_r ,j,i}u^r\ \mod u^q 
\end{equation}
Since we can take $q$ as large as we want, the proof is complete.
\end{pf}
\subsection{The universal abelian cover}\lbl{sec.unab}
The considerations of Section \ref{sec.long} extend readily 
 to the universal abelian cover of $X$. The results essentially
correspond to those obtained by Traldi in \cite{Tra} where the role
of string link is replaced by a choice of link projection. We
will omit proofs since they are identical to the arguments in Section
\ref{sec.long}. In fact the results in Sections \ref{sec.Hom}  and
\ref{sec.Alex} on the infinite cyclic covering and the one-variable
polynomial are consequences of the analogous results for the
universal abelian covering and the multivariable polynomial, but we
will need the details of the argument in the one-variable case for
the arguments in Sections \ref{sec.ks} and after.

 The universal abelian covering
$p:\hat X\to X$ is defined by the epimorphism $\pi\to\ZZ^m$ sending
$\mu_i\to t_i$, where
$\ZZ^m$ is the free abelian (multiplicative) group with generators
$t_1 ,\cdots t_m$. We consider the relative homology group $H_1 (\hat
X,\hat\ast )$, where
$\ast$ is the base-point of $X$ and $\hat\ast =p^{-1}(\ast )$.
This is a module over the Laurent polynomial ring $\zti$. Consider
the multiplicative subset $\S\sub\zti$ consisting of all $f(t_1
,\cdots ,t_m )$ such that $f(1,\cdots ,1))=1$, and the {\em
localization }
$\zti_{\S}$ of
$\zti$ consisting of all quotients of elements of $\zti$ by elements
of
$\S$. We also consider the {\em completion } of
$\zti$. Let $\zvi$ denote the ring of power series in the variable
$v_1 ,\cdots ,v_m$. Then we can define a ring homomorphism $\hat\th
:\zti\to\zvi$ by
$\hat\th (t_i)=1+v_i$. Since $\hat\th (\S )\sub\text{ units of
}\zvi$,
$\hat\th$ extends to an imbedding $\zti_{\S}\to\zvi$ which we still
denote by $\hat\th$.
\begin{lemma}\lbl{lem.freei} $H_1 (\hat X,\hat\ast )_\S$ is a free
$\zti_{\S}$-module with basis $\hat\mu_1 ,\ld ,\hat\mu_m$, where
$\hat\mu_i$ is represented by a lift of $\mu_i$ to a curve in $\hat
X$ starting at some prescribed base-point in $\hat\ast$.
\end{lemma}

We define $\hat\l_i\in H_1 (\hat X,\hat\ast )$ to be the
class represented by the lift of $\l_i$ to a path beginning at
the chosen base point. Then we can write, in $H_1 (\hat X,\hat\ast
)_{\S}$:
\begin{equation}\lbl{eq.ciji}
\hat\l_i =\sum_{j=1}^m \hat c_{ij}^S\hat\mu_j ,\quad 1\le i\le m
\end{equation}
We now have
\begin{lemma}\lbl{lem.ciji} 
\begin{equation}\lbl{eq.thciji}
\hat\th (\hat c_{ij}^S)=\sum_{i_1 ,\ld ,i_k}\mb_{i_1 ,\ld ,i_k
,j,i}(S)v_{i_1}\cdots v_{i_k}
\end{equation}
\end{lemma}

\subsection{Relations in the longitudinal matrix}
We now point out that the matrices $(c_{ij}^S)$ and
$(\hat c_{ij}^S)$ are degenerate. 
\begin{lemma}\lbl{lem.deg} 
$$\sum_{i=1}^m c_{ij}^S=0=\sum_{j=1}^m c_{ij}^S $$
\end{lemma}
\begin{pf} First recall that the $\{\l_i\}$ satisfy the relation which, with our
orientation convention, reads:
\begin{equation}\lbl{eq.lami}
(\l_1^{-1}\m_1^{-1}\l_1 )(\l_2\m_2\l_2^{-1} )(\l_3^{-1}\m_3^{-1}\l_3 )\cdots
=\m_1^{-1}\m_2\m_3^{-1}\cdots
\end{equation}
This is apparent since both sides are represented by $0\times S^1\sub
X$.  If we apply the free differential calculus to this
relation we obtain the following relation in $H_1 (\ti X,\ti\ast )$:
\begin{equation*}
\sum_{i\text{ odd}} t^{-l_i}((t^{-1}-1)\ti\l_i -t^{-1}\ti\m_i
)+\sum_{i\text{ even}}((t^{-1}-1)\ti\l_i
+t^{-1}\ti\mu_i)=\sum_{i=1}^m(-1)^{i}t^{-1}\ti\m_i
\end{equation*}
where $l_i$ is the total linking number of $\l_i$. Since we have
chosen $\l_i$ so that $l_i =0$, this equation becomes 
$(t^{-1}-1)\sum_{i=1}^m \ti\l_i =0$.
Since we are in a free module this becomes, simply,
$\sum_{i=1}^m\ti\l_i =0$, which proves the first equality.

Consider the boundary operator $\pa :H_1 (\ti X,\ti\ast )\to H_0
(\ti\ast )\iso\zt$ from the homology sequence of $(\ti X,\ti\ast
)$. Then $\pa (\ti\m_i )=t-1$ and $\pa (\ti\l_i )=t^{l_i}-1$.
Thus we obtain the equality 
$$t^{l_i}-1=\pa (\ti\l_i )=\sum_j c_{ij}^S\pa (\ti\m_j )=(t-1)\sum_j
c_{ij}^S $$
Since $l_i =0$ and $t-1$ is a non-zero divisor, we obtain the
second inequality.
\end{pf}

The degeneracy for
$(\hat c_{ij}^S)$ is somewhat more complicated, but the argument is
identical to that for Lemma \ref{lem.deg}.
\begin{lemma}\lbl{lem.degi}
\begin{enumerate}
\item $$\tau_i -1=\sum_{j=1}^m (t_j -1)\hat c_{ij}^S$$
where $\tau_i =\prod_j t_j^{l_{ij}}$ and $l_{ij}$ is the linking
number of $\l_i$ and $L_j$.
\item $$a_j (\tau_j -1)=\sum_{i=1}^m a_i (t_i -1)\hat c_{ij}^S$$
where 
$a_i =b_i\cdot\prod\limits_{r\text{
odd}<i}t_r^{-1}\cdot\prod\limits_{r\text{ even}<i}t_r$ 
and 
$b_i =\cases 1& \text{if $i$ even}\\
\tau_i^{-1}t_i^{-1} & \text{if $i$ odd}\endcases$.
\end{enumerate}
\end{lemma}
\begin{pf} Consider the boundary operator $\pa :H_1 (\hat X,\hat\ast
)\to H_0 (\hat\ast )\iso\zti$ from the homology sequence of $(\hat
X,\hat\ast )$. Then $\pa (\hat\m_i )=t_i-1$ and $\pa (\hat\l_i
)=\tau_i -1$ and the first equality follows.

The second equality is proved by applying the Fox differential
calculus to Equation\eqref{eq.lami} as in the proof of Lemma
\ref{lem.deg}.
\end{pf}

\section{First results on the Alexander polynomial}\lbl{sec.Alex}
\subsection{A presentation of the Alexander module}\lbl{sec.al}
Let $Y$ denote the complement of $L_S$, the closure of $S$.
Then $Y=X\cup X_0$, where $X_0$ is the complement of the trivial
$m$-component string link. One sees immediately that $\pi_1
(Y)\iso\pi_1 (X)/<[\m_i ,\l_i ]>$ and, therefore, that:
\begin{equation}\lbl{eq.hom}
H_1 (\ti Y,\ti\ast )\iso H_1 (\ti X ,\ti\ast )/M
\end{equation}
where $M$ is the submodule generated by the elements $(1-t)\ti\l_i
-(1-t^{l_i})\ti\m_i$. Since $l_i =0$ we conclude that $H_1 (\ti
Y,\ti\ast )$ has a presentation with $m$ generators $\{\ti\m_i\}$ and
$m$ relators $\{ (1-t)\sum_j c_{ij}^S\ti\m_j\}$.
\subsection{The Alexander polynomial}
Recall the definition of the Alexander polynomial $\D_L
(t)$ of a link $L$ as a generator of the order ideal of $H_1 (\ti
Y)$, where $Y$ is the complement of $L$. It follows from the
homology exact sequence of $(\ti Y.\ti\ast )$ that $H_1 (\ti
Y,\ti\ast)\iso H_1 (\ti Y)\oplus\zt$ and so the order ideal of $H_1 (\ti
Y)$ is the same as the ideal generated by the $(k\times k)$ minors
of the matrix of a presentation of $H_1 (\ti Y,\ti\ast)$ with
$k+1$ generators. Thus we can conclude that the image of
$\D_{L_S}(t)$ in $\zts$ is a generator of the ideal generated by
the $(m-1)\times (m-1)$ minors of the matrix $(t-1)c_{ij}^S$. Since,
by  Lemma \ref{lem.deg}, the sum of the rows and the sum of the
columns is zero, we have:
\begin{lemma}\lbl{lem.al} The image of $\D_{L_S}(t)$ in $\zts$ is, up to
multiplication by an element of $\S$,
$ (t-1)^{m-1}\det C^S $, where $C^S$ is the $(m-1)\times (m-1)$
matrix whose entries are $\{ c_{ij}^S,\ 1\le i,j\le m-1\}$.
\end{lemma}

Now define $\D_L (u)=\ti\th (\D_L (t))\in\zu$. Then, as a consequence
of Lemmas \ref{lem.cij} and \ref{lem.al}, we have:
\begin{thm}\lbl{th.Phi} There exists a rational power series $\th
(u)$ with constant term $\pm 1$ such that
$$ \D_{L_S}(u)=\Phi_S (u)\th (u) $$
where $\Phi_S (u)$ is defined by:
\begin{equation}\lbl{eq.Phi}
\Phi_S (u)=u^{m-1}\det (\sum_{1\le i_1 ,\ldots ,i_r\le m}\mb_{i_1
,\ldots ,\i_r ,j,i}(S)u^r )
\end{equation}
\end{thm}
A {\em rational } power series is one which is the expansion of a
quotient of polynomials in $u$. For example, the elements of
$\im\ti\th$ are rational. See \cite{Tra1} for a result about the
multivariable Conway polynomial which is similar to Theorem \ref{th.Phi}.

From this theorem we can see, in particular, that the first
non-vanishing coefficient of $\D_{L_S}(u)$ is, up to sign, equal
to the first non-vanishing coefficient of $\Phi_S (u)$. In
particular we have proved Corollary \ref{cor.lnk}, for
$\D_{L_S}(u)$ rather than the Conway polynomial, up to sign.

We can also obtain Corollary \ref{cor.conc} for $\D_{L_S}(u)$, up
to sign\lbl{conc}. Suppose $L_1$ and $L_2$ are concordant. We may
assume, by a theorem of Tristram \cite{Tri}, that $L_2$ is obtained
from
$L_1$ by a ribbon move. In other words, for some trivial link $T$ in a
ball $B$ disjoint from $L_1$, we obtain $L_2$ by band-summing each
component of $T$ to some component of $L_1$. Now we can lift $L_1$ to a
string link $S_1$ by choosing a so-called {\em d-base } (see Habegger-Lin
\cite{HaL}), i.e. an imbedded $2$-disk which meets each component of $L_1$ in
a single point. Clearly we can choose a d-base which is disjoint from
$B$ and the bands used to obtain $L_2$. So the same
d-base can be used for $L_2$ and lifts $L_2$ to a string link $S_2$ concordant
to
$S_1$. Since the $\mb$-invariants of a string link are concordance
invariants (see \cite{Sta}), $\Phi_{S_1}(u)=\Phi_{S_2}(u)$. Thus it
follows from Theorem \ref{th.Phi} that the first non-vanishing
coefficients of $\D_{L_i}(u)$ are the same (up to sign). Finally,
connect sum of a link with a local knot clearly multiplies the Alexander
polynomial of the link by the Alexander polynomial of the knot. 
\subsection{The multivariable Alexander polynomial}\lbl{sec.mual}
The same argument as in Section \ref{sec.al} shows that $H_1 (\hat
Y ,\hat\ast )$ has a presentation with generators $\hat\mu_1
,\cdots ,\hat\mu_m$ and relations 
$$(\tau_i -1)\hat\mu_i =(t_i -1)\hat\l_i =\sum_{j=1}^m (t_i -1)\hat
c_{ij}^S\hat\mu_j $$

The usual definition of the multivariable Alexander polynomial
$\D_L (t_1 ,\cdots ,t_m )$ is the greatest common divisor of the
ideal $E$ generated by the $(k\times k)$-minors of a
presentation matrix with $(k+1)$ generators of the $\zti$-module
$H_1 (\hat Y ,\hat\ast )$. We have the presentation matrix 
\begin{equation}\lbl{eq.pres}
P=((t_i -1)\hat c_{ij}^S -\delta_{ij}(\tau_i
-1))
\end{equation}
 for $H_1 (\hat Y ,\hat\ast )_{\S}$ and so $E_{\S}$ is generated
by the $(m-1)\times (m-1)$-minors of $P$. By Lemma
\ref{lem.degi} the sum of unit multiples of the rows and the sum of
unit multiples of the columns is zero, we conclude that $E_{\S}$
is principal generated by any $(m-1)\times (m-1)$-minor $\D$ of
$P$. Since $\zti_{\S}$ is a localization of $\zti$, it follows
that $\D_{L_S}(t_1 ,\cdots t_m )$ is a unit multiple of $\D$ in
$\zti_{\S}$. Theorem \ref{th.mual} now follows from \eqref{eq.pres} and
Lemma \ref{lem.ciji}.

The remainder of this paper will focus on the one-variable
polynomial and, in particular, be devoted to sharpening Theorem
\ref{th.Phi} to obtain Theorem \ref{th.slnk}

\section{The Conway polynomial}\lbl{sec.ks}
We recall the definition of the Conway polynomial of a link $L$ in
$\R^3$ (see, for example, \cite{Lev}). Let $V$ be a Seifert surface for
$L$, i.e.
$V\sub\R^3$ is an oriented surface with $\pa V =L$. We define the
Seifert pairing
$\s :H_1 (V)\times H_1 (V)\to\Z$ by $\s (\a ,\b )=\l (j_{\ast}\a
,\b )$ where $j:V \to\R^3 -V$ is defined by a push in the
positive normal direction and $\l :H_1 (\R^3 -V )\times H_1 (V
)\to\Z$ is the linking pairing, which is non-singular by
Alexander duality. Let $A$ be a matrix representing
$\s$ with respect to a basis of $H_1 (V)$. Then define
the {\em potential function } $\O_L (t)=\det (tA-t^{-1}A^t
)$. This Laurent polynomial depends only on $L$ and  $\O_L (t)$ is a
polynomial in $t-t^{-1}$. Thus we may define the {\em Conway polynomial
} $\nl$ by the equation $\n_L (t-t^{-1})=\O_L (t)$.

Now suppose we are given a string link $S=\{ S_i\}$. We have already discussed 
the link closure $L_S$ in Section \ref{sec.str}; now we associate an oriented
knot
$K_S$ to $S$ by closing it in the following different manner, as indicated in
Figure \ref{fig.2}.

\begin{figure}[hbt]
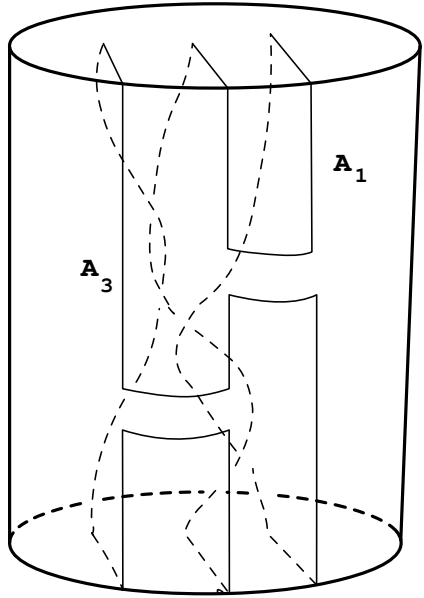

\begin{center}{\BoxedEPSF{Fig2 scaled 600}}\caption{The knot
closure of a string link}\lbl{fig.2}
\end{center}
\end{figure} 
Insert bands $B_i\approx I\times I$, for $1\le i\le m-1$, into
$I\times S^1$ so that $B_i$ connects $S_i$ to $S_{i+1}$ and
$B_i\cap S_i\ap I\times 0$ and $B_i\cap S_{i+1}\ap I\times 1$. We
also arrange them so that $B_{i+1}$ is below $B_i$. Then $K_S$
is obtained from $L_S\cup\bigcup_i B_i$ by removing the part of each
$B_i$ corresponding to $(0,1)\times I$ and orienting it consistent
with the orientation of $L_S$. Now choose a Seifert surface $W$
for $K_S$ in $I\times D^2$ so that $W\cap\pa (I\times D^2 )=K_S\cap\pa
(I\times D^2 )$. By the simple modification of merely adjoining the bands $\{
B_i\}$ to $W$, we can convert $W$ into a Seifert surface $V$ for $L_S$ (see
Figure \ref{fig.3}).

\begin{figure}[hbt]
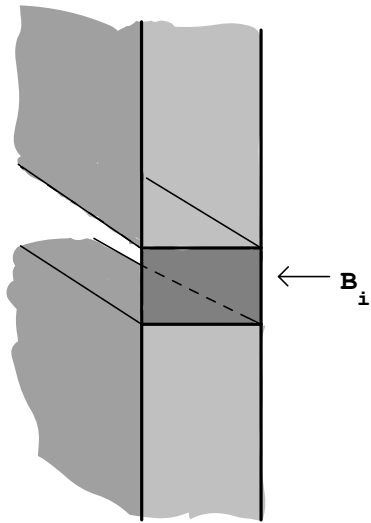

\begin{center}{\BoxedEPSF{Fig3 scaled 600}}\caption{Adjoining a
band to a Seifert surface}\lbl{fig.3}
\end{center}
\end{figure} 
The Seifert pairings and
matrices of $W$ and $V$ are closely related. It is clear that the
inclusion $W\sub V$ induces an isomorphism:
\begin{equation}\lbl{eq.VW}
H_1 (V)\iso H_1 (W)\oplus\Z^{m-1}
\end{equation}
The $\Z^{m-1}$ summand has, as basis, the classes represented by any
$m-1$ of the components $\{ L_i\}$ of $L_S$. We will usually
choose $L_1 ,\ldots ,L_{m-1}$. The Seifert pairing of $V$, when
restricted to $W$, obviously coincides with that of $W$. Therefore
if $A$ is the Seifert matrix of $W$, with respect to some basis of
$H_1 (W)$, then the Seifert matrix $\ba$ of $V$ will be of the
form:
\begin{equation}\lbl{eq.sm}
\ba =\pmatrix A & B^t \\ B & \L \endpmatrix
\end{equation}
where $M^t$ denotes the transpose of the matrix $M$ and $\L
=(l_{ij}),\ 1\le i,j\le m-1$ where (for $1\le i,j\le m$) 
\begin{equation}\lbl{eq.lij}
 l_{ij}=\cases \text{linking number of }L_i \text{ and
}L_j &
\text{ if } i\not= j \\ -\sum_{r\not= i}l_{ir} & \text{if } i=j \endcases 
\end{equation} 
Note that the parallel push-off of each $L_i$ has total linking
number $0$. The Conway polynomials of $K_S$ and $L_S$ are given
by:
\begin{equation}\lbl{eq.Con}
\begin{split} \n_{K_S}(s-s^{-1}) & =\det (sA-s^{-1}A^t ) \\
\n_{L_S}(s-s^{-1}) & =\det (s\ba-s^{-1}\ba^t )
\end{split}
\end{equation}
Since $E=A-A^t$ is unimodular, it follows that $\n_{K_S}(0)=1$, as
usual. Setting $z=s-s^{-1}$, the Conway polynomials lie in the
polynomial ring $\zz\sub\zs$. Let $\S '$ denote the multiplicative
subset of $\zs$ consisting of all $f(s)$ satisfying $f(1)=\pm
f(-1)=\pm 1$. Then, over the localized ring $\zss$, the matrix
$\SS =sA-s^{-1}A^t =zA-s^{-1}E$ is invertible. We can, therefore, by
elementary row operations in $\zss$ convert $s\ba -s^{-1}\ba^t$
into
\begin{equation}\lbl{eq.sbar}
\bar{\SS} =\pmatrix \SS & B^t \\ 0 & z\L -z^2 B\SS^{-1}B^t
\endpmatrix 
\end{equation}
From this we obtain
\begin{equation}\lbl{eq.fact}
\n_{L_S}(z)=\n_{K_S}(z)\Ga_S (z)
\end{equation}
where 
\begin{equation}\lbl{eq.phi}
\Ga_S (z)=\det (z\L -z^2 B\SS^{-1}B^t )
\end{equation}
Note that, although $\Ga_S$ is, by its definition, an element of
$\zss$, it follows from \eqref{eq.fact} that it is a quotient of
polynomials in $z$ or, alternatively, a power series in $z$.
\section{The Seifert matrix and the $\mb$-invariants}\lbl{sec.Seif}
\subsection{Homology from the Seifert pairing}\lbl{sec.seif}
We would like to relate the matrices $\SS$ and $\bar{\SS}$ to
the longitudinal matrix $(c_{ij}^S)$ defined in Equation
\eqref{eq.cij}. For that we have to describe the procedure by
which one defines a presentation for $H_1 (\ti X,\ti\ast )$,
using a Seifert matrix. This will be entirely analogous to the
traditional way of using a Seifert matrix of a knot or link $K$
to produce a presentation for $H_1 (\widetilde{S^3 -K})$ (see \cite{Kau}). 

Let $S, W, V$ be as above and let $Y=I\times D^2-V$. Let $i_+
,i_- :V\to Y$ be defined by a push-off in the positive or
negative normal direction, respectively. Now $V$ lifts into
$\ti X$ and the translates $\{ V_i\}$ of $V$ cut $\ti X$
into the union of the translates $\{ Y_i\}$ of a lift of $Y$.
We may WLOG move the base-point $\ast$ slightly into the
interior of $Y$ and assume $\ast\in V$. Let $B$ be a small
ball containing $\ast$ such that $B\cap V=B'$ and $\pa B\cap
Y=B_+\cup B_-$, where $B', B_+ ,B_-$ are $2$-disks. Now $\ti
X-\int{\ti B}$ is a union of the $\{ Y_i\}$  (now redefined as
the lifts of $Y-(\int B\cap Y)$) attached along $\{ V_i
-B_i\}$, where the $\{ B_i\}$ are the lifts of $B'$.

A standard Mayer-Vietoris argument produces an exact sequence of\hfil\break
$\zt$-modules:
\begin{equation}\lbl{eq.MV}
\zt\otimes H_1 (V-B')\stackrel{\rho}{\rightarrow}\zt\otimes H_1 (Y,B_+\cup B_-
)\to  H_1 (\ti X,\ti B)\to 0   
\end{equation}
where $\rho$ is the $\zt$-homomorphism defined by:
$$\rho (1\otimes\a )=t\otimes i_{+\ast}(\a )-1\otimes i_{-\ast}(\a
)$$

Alexander duality gives us an isomorphism  
$$H_1 (Y,B_+\cup B_- )\iso H^1 (V-B',V\cap\pa (I\times D^2 ))$$ 
which is adjoint to the non-singular pairing:;
$$ \LL :H_1 (Y,B_+\cup B_- )\times H_1 (V-B',V\cap\pa (I\times D^2 ))\to\Z
$$
defined by linking number. Paths representing elements of these homology
groups can be closed up (disjointly) in an obvious way.

From $\LL$ we can define two Seifert pairings:
$$ \s_{\pm}: H_1 (V-B')\times H_1 (V-B',V\cap\pa (I\times D^2 ))\to\Z $$
by $\s_{\pm}(\a ,\b )=\LL (i_{\pm\ast}(\a ),\b )$. To get Seifert
matrices we need to choose bases for $H_1 (V-B')$ and $H_1 (V-B',V\cap\pa
(I\times D^2 ))$. The basis for $H_1 (V-B')$  can be the same basis we used for
$H_1 (V)$ in \eqref{eq.sm} with the addition of the class of $L_m$. For $H_1 (V-B',V\cap\pa
(I\times D^2 ))$ we note that $V\cap\pa(I\times D^2 )$ consists of $m$ disjoint
arcs, one on each component of $\pa V=L_S$. Thus $H_1 (V-B',V\cap\pa (I\times D^2
))\iso H_1 (V-B')\oplus\Z^{m-1}$. For basis we can choose the one already chosen
for $H_1 (V-B')$ together with the classes of the arcs $a_1 ,\cdots ,a_{m-1}$,
where $a_i$ crosses the band $B_i$ connecting $L_i\cap\pa (I\times D^2 )$ to
$L_{i+1}\cap\pa (I\times D^2 )$ (see Figure \ref{fig.4}).
 
\begin{figure}[hbt]
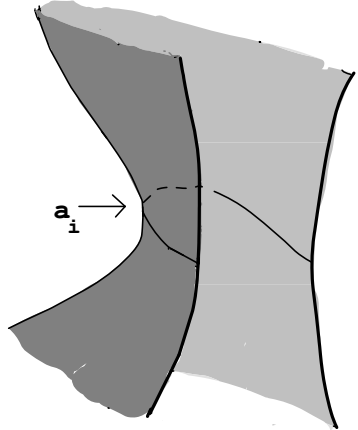

\begin{center}{\BoxedEPSF{Fig4 scaled 600}}\caption{Definition
of $a_i$}\lbl{fig.4}
\end{center}
\end{figure}
Let $A_{\pm}$ denote the matrix representing $\s_{\pm}$ with respect to these
bases. So $tA_+ -A_-$ (with entries in $\zt$) represents the map $\rho$ in
\eqref{eq.MV} and so is a presentation matrix of $H_1 (\ti X,\ti\ast )$. Note
that $A_{\pm}$ has the form  $A_+ =\pmatrix
\ti A  & U_+
\endpmatrix , A_- =\pmatrix \ti A^t  & U_- 
\endpmatrix$, where $\ti A$ is the matrix determined by  $\bar A$
from
\eqref{eq.sm} as follows:
\begin{equation}
\ti A = \pmatrix A & \ti B^t \\ \ti B & \L \endpmatrix
\end{equation}
and $\ti B$ is the matrix determined from $B$ by:

\begin{equation}\lbl{eq.tiA}
\begin{split}
(i) & \text{  The submatrix obtained by deleting the last row is $B$. }\\
(ii) & \text{ The sum of the rows of $\ti B$ is $0$. }
\end{split}
\end{equation}
and $\L$ is the matrix with entries $l_{ij}, 1\le i,j\le m$ (see
\eqref{eq.lij}). Assertion (ii) follows because $\sum_i L_i$ bounds $V$.
 $U_{\pm}$ is a matrix with $m-1$ columns whose entries are determined by
 the following values of the Seifert pairings:
\begin{equation}\lbl{eq.seif}
\begin{split}
\s_{\pm}(\a ,a_j ) &=0\quad \text{ if }\a\in H_1 (W) \\
t_{ij}^+ =\s_+ (L_i ,a_j ) &= \cases 0\quad\text{ if $j$ is even }  \\
1 \quad\text{ if $i=j$ is odd }  \\
-1 \quad\text{ if $i-1=j$ is odd }\endcases \\
t_{ij}^- =\s_- (L_i ,a_j ) &=\cases 0\quad\text{ if $j$ is odd }  \\
1 \quad\text{ if $i=j$ is even }  \\
-1 \quad\text{ if $i-1=j$ is even }\endcases \\
\end{split}
\end{equation}
These assertions follow from the fact that $a_j$ bounds a small disk
disjoint from
$W$ and any  $L_i$, for $i\not= j,j+1$ (see Figure \ref{fig.5}).

\begin{figure}[hbt]
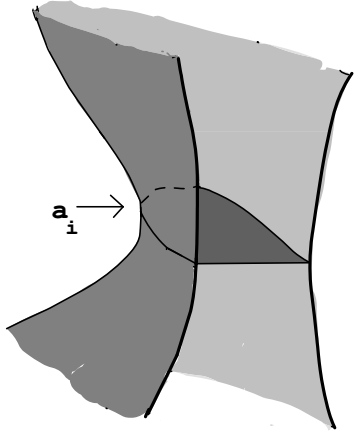

\begin{center}{\BoxedEPSF{Fig5 scaled 600}}\caption{The disk
bounded by $a_i$}\lbl{fig.5}
\end{center}
\end{figure}
Then we can write $U_{\pm}= \pmatrix 0  \\ T_{\pm} \endpmatrix$ where $T_{\pm}=
(t_{ij}^{\pm})$ is the
 $m\times (m-1)$-matrix defined in \eqref{eq.seif}.

Now the matrix $\SS_t =tA-A^t$ is invertible in $\zts$ and so, by elementary
row operations, we can convert the matrix $tA_+ -A_-$ to:
\begin{equation}\lbl{eq.st}
\hat{\SS_t} = \pmatrix \SS_t & (t-1)\ti B^t & 0 \\ 0 & (t-1)\L  -(t-1)^2\ti B
\SS_t^{-1}\ti B^t & tT_+ -T_- \endpmatrix 
\end{equation}
 $\hat{\SS_t}$ is now a presentation matrix for $H_1 (\ti X,\ti\ast
)_{\S}$ whose generators are the duals in $H_1 (Y,B_+\cup B_- )$, with respect
to the pairing $\LL$, of
the given basis $\{\b_i , L_i , a_i\}$ of $H_1 (V-B',V\cap\pa (I\times D^2 ))$
($\b_i$ denotes a basis of $H_1 (W)$). Let us consider the elements
$\ti\l_i , \mu_i '\in H_1 (Y,B_+\cup B_- )$, where $\ti\l_i
=i_{+\ast}(L_i)$, as in \eqref{eq.cij}, and $\mu_i '$ is the meridian
curve which starts at $B_+$ and travels along a positive push-off of a
curve $\gg_i$ in $W$ to $L_i$, half-way round a small meridian of
$L_i$ and back to $B_-$ along a negative push-off of $\gg_i$
(see Figure \ref{fig.6}). 

\begin{figure}[hbt]
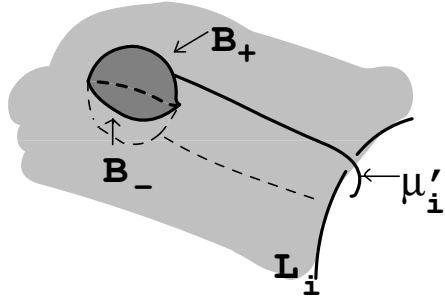

\begin{center}{\BoxedEPSF{Fig6 scaled 600}}\caption{The meridian
curve $\mu_i'$}\lbl{fig.6}
\end{center}
\end{figure}
Then we have the following values of $\LL$:
\begin{equation}\lbl{eq.LL}
\begin{split}
\LL (\ti\l_i ,L_j ) &= l_{ij} \\
\LL (\ti\l_i ,a_j ) &= t_{ij}^+ \\
\LL (\ti\l_i ,\b_j ) &=\ti b_{ij} \\
\LL (\mu_i ' ,L_j ) &= \delta_{ij}  \\
\LL (\mu_i ' ,a_j ) &=0 \\
\LL (\mu_i ' ,\b_j ) &=0 
\end{split}
\end{equation}
The last lines follows from the fact that $\gg_i$ lies in $W$.

From \eqref{eq.LL} we can deduce $L_i^{\sharp} = \mu_i '$ and:
\begin{equation}\lbl{eq.dual}
\ti\l_i = \sum_j t_{ij}^+a_j^{\sharp} +\sum_jl_{ij}L_j^{\sh}+\sum_j\ti
b_{ij}\b_j^{\sh}
\end{equation}
From the presentation matrix \eqref{eq.st} we can write:
\begin{equation}
\begin{split}
\sum_j s_{ij}\b_j^{\sh} + (t-1)\sum_j \ti b_{ji}L_j^{\sh} & =0 \\
(t-1)\sum_j (l_{ij}-(t-1)\th_{ij})L_j^{\sh}  +\sum_j\tau_{ij}a_j^{\sh} &=0
\end{split}
\end{equation}
where $\SS_t =(s_{ij}),\ \ti B\SS_t^{-1}\ti B^t =(\th_{ij})$ and $tT_+ -T_-
=(\tau_{ij})$.

Note that $\SS_t$ is invertible and we can see from \eqref{eq.seif} that the
$(m-1)\times (m-1)$ submatrix $\TT$ obtained by deleting the last row of
$tT_+ -T_-$ is invertible. Thus we may write:
\begin{equation}\lbl{eq.sharp}
\begin{split}
\b_i^{\sh} &= (1-t)\sum_{r,j}\bar s_{ir}\ti b_{jr}L_j^{\sh} \\
a_i^{\sh} &= (1-t)\sum_{r,j}\bar\tau_{ir}(l_{rj}-(t-1)\th_{rj})L_j^{\sh}
\end{split}
\end{equation}
where $\SS_t^{-1}=(\bar s_{ij})$ and $\TT^{-1}=(\bar\tau_{ij})$ for $j\le
m-1$. For convenience later we define $\bar\tau_{im}=0$ and let $\ti\TT_t
=(\bar\tau_{ij})$, for $j\le m$.

Substituting \eqref{eq.sharp} into \eqref{eq.dual} we obtain:
\begin{equation*}
\ti\l_i =\sum_j (
(1-t)\sum_{s,r}t_{is}^+\bar\tau_{sr}(l_{rj}-(t-1)\th_{rj})+l_{ij}+(1-t)\ti
b_{is}\bar s_{sr}\ti b_{jr})L_j^{\sh}
\end{equation*}
Substituting the definition of $\th_{ij}$ this becomes:
\begin{equation}
\ti\l_i =\sum_j (
(1-t)\sum_{s,r}t_{is}^+\bar\tau_{sr}(l_{rj}-(t-1)\th_{rj})
+l_{ij}-(t-1)\th_{ij}) L_j^{\sh}
\end{equation}
or, a bit more succinctly:
\begin{equation}\lbl{eq.lmi}
\ti\l_i =\sum_j ( \sum_r (\d_{ir}-(t-1)\sum_s
t_{is}^+\bar\tau_{sr})(l_{rj}-(t-1)\th_{rj})) L_j^{\sh}
\end{equation}
We can rewrite \eqref{eq.lmi} briefly as:
\begin{equation}\lbl{eq.Cbij}
 \ti\l_i =\sum_j C_{ij}\mu_j ' 
\end{equation}
where we define the matrix $C$ by: 
\begin{equation}\lbl{eq.CC}
C=(I-(t-1)T_+\ti\TT_t )(\L -(t-1)\ti B\SS_t^{-1}\ti B^t )
\end{equation}

We now need to compare $C$ to $(c_{ij}^S )$ from \eqref{eq.cij}. To do
this we must compare the two sets of meridian generators $\{\ti\mu_i\}, \{\mu_i
'\}$ of $H_1 (\ti X,\ti\ast )$.  The difference between the representative
closed curves, $\{\ti u_i\}$ and $ \{ u_i
'\}$, respectively, in $X$ is that the stems of $u_i$ lie in $0\times
D^2$ while the stems of $u_i '$ lie essentially in $W$ (see Figure \ref{fig.7}).

\begin{figure}[hbt]
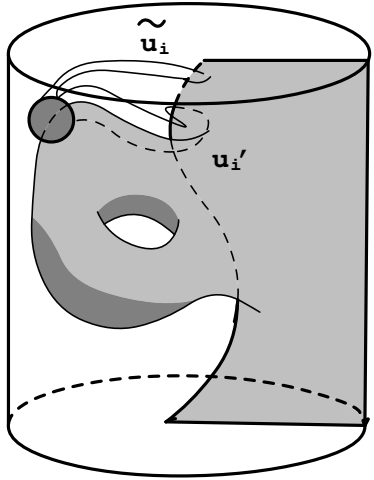

\begin{center}{\BoxedEPSF{Fig7 scaled 600}}\caption{The two
choices of meridian curves}\lbl{fig.7}
\end{center}
\end{figure}
However if we choose $W$ appropriately, then we can also choose $u_i '$ to be a
slight translate of $u_i$, since we have allowed ourselves to choose the
stems of $u_i '$ arbitrarily in $W$ (see Figure \ref{fig.8}).

\begin{figure}[hbt]
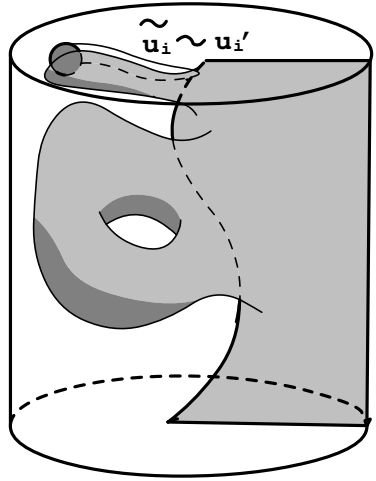

\begin{center}{\BoxedEPSF{Fig8 scaled 600}}\caption{How to make
the two choices of meridian curves homotopic}\lbl{fig.8}
\end{center}
\end{figure}
	Thus we will have succeeded in showing that $C=(c_{ij}^S )$ (with the
substitution $t=1+u$) if we prove:
\begin{lemma}\lbl{lem.seif}
The matrix $\ti B\SS_t^{-1}\ti B^t$ depends only on $S$, i.e. it is
independent of the choice of $W$.
\end{lemma}
\begin{pf} Any two choices of $W$ are cobordant so it follows from the usual
argument (see e.g. \cite{Lev1}) that the Seifert matrices are S-equivalent. More
precisely there are a sequence of elementary cobordisms connecting the two
choices of $W$, where an elementary cobordism consists of adjoining or
removing the boundary of a solid torus. The addition produces two new
generators $\xi ,\eta$ which are, respectively, the meridian and longitude
curves of the torus.   This changes $A$ by the following enlargement. 
\begin{equation}\lbl{eq.A}
 A\longrightarrow \pmatrix A & \vdots & 0 \\ \hdots &
\cdot & 0
\\ 0 & 1 & 0
\endpmatrix  
\end{equation}
where the last row and column  corresponds to $\xi$. The effect on
$\ti B$ is the following enlargement.
\begin{equation}\lbl{eq.B}
 \ti B\longrightarrow\pmatrix \ti B & \cdot & 0
\endpmatrix
\end{equation}
 From \eqref{eq.A} we see that the resulting enlargement of $\SS_t
=tA-A^t$ is:
$$ \SS_t \longrightarrow \pmatrix \SS_t & \vdots &0 \\ \hdots & \cdot & -1 \\
0 & t & 0 \endpmatrix $$
From this we see that the resulting enlargement of $\SS_t^{-1}$ is:
\begin{equation}\lbl{eq.S-1} \SS_t^{-1}\longrightarrow \pmatrix \SS_t^{-1} & 0 &
\vdots
\\ 0 & 0 & t^{-1}  
\\ \hdots & -1 & \cdot \endpmatrix 
\end{equation}
Combining \eqref{eq.B} and \eqref{eq.S-1} and carrying out the matrix
multiplication we can check that $\ti B\SS_t^{-1}\ti B^t$ is unchanged by these
enlargements of $\SS_t^{-1}$ and $\ti B$. 

This completes the proof of Lemma \ref{lem.seif}.
\end{pf} 
\subsection{Completion of proof of Theorem \ref{th.slnk}}
We have now laid the necessary groundwork to relate the matrix   
$\Th (z)= z\L -z^2  B\SS^{-1}B^t$ from \eqref{eq.sbar},  
whose determinant is $\Ga_S (z)$ (see
\eqref{eq.phi}) and the matrix $(c_{ij}^S )$ from \eqref{eq.cij} which is a
power series in $u=t-1$ whose coefficients are given, in Lemma \ref{lem.cij},
in terms of the
$\mb$-invariants of
$S$. Define $\ti\Th (z)=z\L -z^2\ti B\SS^{-1}\ti B^t$, so that
$\Th (z)$ is obtained from $\ti\Th (z)$ by deleting the last row and column,
where $z=t-t^{-1}$. Then, using equation \eqref{eq.CC}, we have:
\begin{equation}
\begin{split}
\ti\Th (z) &=z(\L -(t^2 -1)\ti B\SS_{t^2}^{-1}\ti B^t )\\
&=z(I-(t^2 -1)T_+\ti\TT_{t^2})^{-1}(c_{ij}^S (u))
\end{split}
\end{equation}
where we are now using the substitution $u=t^2 -1$. We now need to point out
that $I-(t-1)T_+\ti\TT_t$ is a lower triangular matrix whose diagonal entries
are $t^{-1},1,t^{-1},1,\ldots$. It is left to the reader to verify this from the
definitions- see \eqref{eq.seif}. From this we can now write:
\begin{equation*}
\det\Th (z)= u^{m-1}t^{\epsilon}\det (c_{ij}^S )_{1\le i,j\le m-1} 
\end{equation*}
where $\epsilon =\cases 0\quad\text{if $m$  is odd}\\ 1\quad\text{if $m$ is
even}\endcases$ and $u=t^2 -1=tz$. Substituting $t=\sqrt{u+1}$ and
$z=u/\sqrt{u+1}$, this completes the proof of Theorem
\ref{th.slnk}.
\subsection{Proof of Corollary \ref{cor.con}}
The proof begins with the argument on page~\pageref{conc} in the proof of
Corollary \ref{cor.conc}. Thus we may assume that there are concordant string
links $S, S'$ whose closures are $L, L'$, respectively. It follows
immediately from their definitions, that $K_S$ and $K_{S'}$ are concordant
and, by the classical result of Fox-Milnor \cite{FoM}, there exist $f(z), g(z)$,
as in the statement of Corollary
\ref{cor.con}, such that $\n_{K_S}(z)f(z)f(-z)=\n_{K_{S'}}(z)g(z)g(-z)$. Since
$\Ga_S (z)=\Ga_{S'}(z)$, the  Corollary follows from Theorem \ref{th.slnk}.

\section{Proof of Proposition \ref{prop.gam}}\lbl{sec.prf}
From Equation \eqref{eq.phi} we    have that  $\Ga_S 
(z)=(-1)^{m-1}z^{2(m-1)}\Phi (z)$ where $\Phi (z)=\det A(t)$ and
$A(t)$ is the {\em skew-Hermitian} matrix $B\SS^{-1}B^t$. In this context
$A(t)$ has entries in the localized ring $\zts$ and skew-Hermitian means
the equality $A(t^{-1})=-A(t)^t$. Of course $z=t-t^{-1}$ as usual. It
follows immediately from the skew-Hermitian property that $\Phi (z)
=(-1)^{m-1}\Phi (-z)$ which proves Assertion (1) in Proposition
\ref{prop.gam}.

To prove (2) we place ourselves in the somewhat larger ring $\L =\Q
[t,t^{-1}]_I$, the localization at the principal ideal $I=(t-1)$. 
$\L$ is a {\em discrete valuation ring}, i.e. all its ideals are powers
of $I$ (or rather the ideal generated by $I$ in $\L$). For such a
ring it is a standard argument to show that the skew-Hermitian matrix
$A(t)$ is congruent to a matrix $D(t)$ which is a block sum of
$1\times 1$-matrices  $(\phi (t))$ and $2\times 2$-matrices of the form
$(\smallmatrix 0 & \psi (t) \\ -\psi (t^{-1}) & 0 \endsmallmatrix )$. It
follows that  $\phi (t^{-1})=-\phi (t)$, and so we can write $\phi
(t)=(t-t^{-1})\phi '(t)$, where $\phi '(t^{-1})=\phi '(t)$. Now we have
$\Phi (t)=h(t)h(t^{-1})\det D(t)$. If $D(t)$ has any block summands
$(\phi (t))$, then $\Phi (t)$ is divisible by $t-t^{-1}$. If $m-1$
is even, then we must have another such summand and so $\Phi (t)$ is
divisible by $(t-t^{-1})^2$. In other words $\Phi (z)$ is divisible by
$z^2$. If there are no such summands then it follows that $\Phi (t)$
must be of the form $g(t)g(t^{-1})$. Since all the elements of $\L$ are
rational functions the proof is complete.
  
\section{Example}\lbl{sec.ex} To illustrate Theorem \ref{th.slnk} we consider
a very simple example. Let $S$ be the string link in Figure \ref{fig.9}.

\begin{figure}[hbt]
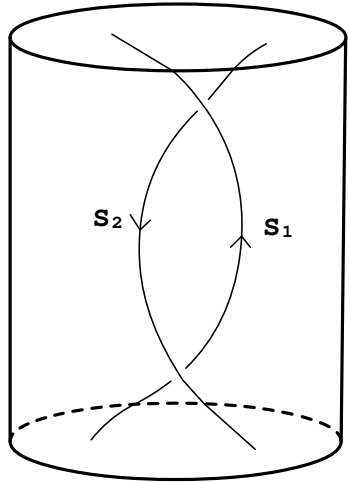

\begin{center}{\BoxedEPSF{Fig9 scaled 600}}\caption{A simple
example}\lbl{fig.9}
\end{center}
\end{figure}
An easy computation gives:
\begin{equation}\lbl{eq.hopf}
\l_1 =\mu_1^{-1}\mu_2\quad\l_2 =\mu_2^{-1}\mu_1
\end{equation}
The $\mb$-invariants we need are computed easily from \eqref{eq.hopf}
and are given by:
$$\mu_{i_1 ,\cdots ,\i_r ,1,1}=\cases (-1)^{r+1}\quad\text{if }i_1 =\cdots
=i_r =1 \\ 0\quad\text{otherwise} \endcases $$
Thus we see that $\l_{11}(u)=-u/u+1$, where $\l_{ij}(u)$ is defined in
Theorem \ref{th.slnk}. From the definition in Theorem \ref{th.slnk} we
see that 
$$\Ga_S (z)=\sqrt{u+1}(-u/u+1)=-u/\sqrt{u+1}=-z $$
Since $K_S$ is trivial, Theorem \ref{th.slnk} tells us that
$\n_{L_S}(z)=-z$, which checks out since $L_S$ is just the right-hand
Hopf link. 

\end{document}